\documentclass[11pt]{amsart}
\usepackage{amsmath,amsfonts,amssymb,amsthm,epsfig,color,tikz}
\newcommand{\bb}{\mathbb}
\newcommand{\p}{\overrightarrow{p}}

\newcommand{\C}{\bb C}

\newcommand{\Z}{\bb Z}
\newcommand{\R}{\bb R}

\newcommand{\Q}{\bb Q}

\newcommand{\G}{\bb G}
\newcommand{\s}{\bb S}

\newcommand{\om}{\omega}
\newcommand{\ww}{\mathbf w}
\newcommand{\Om}{\Omega}

\newtheorem*{ques}{Question}

\newcommand{\F}{\mathcal F}

\newcommand{\vv}{\mathbf{v}}
\newcommand{\uu}{\mathbf{u}}

\newcommand{\minuszero}{\backslash\{0\}}
\newcommand{\La}{\Lambda}
\newtheorem{Theorem}{Theorem}
\numberwithin{Theorem}{section}
\newtheorem{Cor}[Theorem]{Corollary}

\newtheorem{lemma}[Theorem]{Lemma}

\newtheorem*{lemma*}{Lemma}
\newtheorem*{question*}{Question}

\newtheorem*{obs}{Observation}

\newtheorem*{theorem*}{Theorem}

\numberwithin{equation}{section}
\begin{document}
\title{Gap distributions and homogeneous dynamics}
\author{Jayadev S.~Athreya}

\email{jathreya@iilinois.edu}
\address{Deptartment of Mathematics, University of Illinois Urbana-Champaign, 1409 W. Green Street, Urbana, IL 61801}
    \thanks{J.S.A. partially supported by NSF grant DMS 1069153}


\begin{abstract}
We survey the use of dynamics of $SL(2, \R)$-actions to understand gap distributions for various sequences of subsets of $[0, 1)$, particularly those arising from special trajectories of various two-dimensional dynamical systems. We state and prove an abstract theorem that gives a unified explanation for some of the examples we present.\end{abstract}
\maketitle

\section{Introduction}\label{sec:intro} The study of the distribution of gaps in sequences is a subject that arises in many different contexts and has connections with many different areas of mathematics, including number theory, probability theory, and spectral analysis. In this paper, we study gap distributions from the perspective of dynamics and geometry, exploring examples connected with the dynamics of $SL(2, \R)$-actions on moduli spaces of geometric objects, in particular the space of lattices and the space of translation surfaces. 

The inspiration for this article is the a quote from the beautiful paper of Elkies-McMullen~\cite{Elkies}, referring to their explicit computation of the gap distribution of the sequence of fractional parts of $\sqrt{n}$, using the dynamics of the $SL(2, \R)$-action on the space of \emph{affine unimodular lattices} in $\R^2$.
\medskip
\begin{center}\textit{\small ``\ldots the uniform distribution of lattices explains the exotic
distribution of gaps."}\end{center}
\medskip

\noindent Indeed, the main results of our paper, Theorem~\ref{theorem:meta}, Theorem~\ref{theorem:meta:period}, and Theorem~\ref{theorem:meta:circle}, give unified explanations of several examples of `exotic' gap distributions via uniform distribution on various moduli spaces of geometric objects. 
\subsection{Equidistribution, randomness, and gap distributions}\label{subsec:motivate} Suppose that for each positive integer $k$, we are given a finite list of points $F(k) \subset [0, 1)$, where by a {\it list}, we mean a finite
non-decreasing sequence of real numbers where $N_k$ denotes the number of terms in the $k^{th}$ sequence $F(k)$. We write $$F(k) = \left \{F_k^{(0)} \le F_k^{(1)} \le \ldots F_k^{(N_k)}\right \},$$  
and we assume $N_k \rightarrow \infty$ as $k \rightarrow \infty$. In many situations, we are interested in the `randomness' of the sequence of lists $\{F(k)\}_{k=1}^{\infty}$. A first test of `randomness' is whether the lists $F(k)$ \emph{uniformly distribute} in $[0, 1)$, that is the measures $\Delta_k = \frac{1}{N_k} \sum_{j=0}^{N_k} \delta_{F_k^{(j)}}$ converge weak-* to Lebesgue measure, i.e., for any $0 \le a \le b \le 1$,
\begin{equation}\label{eq:equidist}
\lim_{k \rightarrow \infty} \Delta_k(a,b) = b-a.
\end{equation}
A more refined question (not necessarily dependent on (\ref{eq:equidist})) is to examine the distribution of \emph{gaps} for the sequences $F(k)$. That is, form the associated \emph{normalized gap sets}
\begin{equation}\label{eq:gaps}
G(k) : = \left\{ N_k\left(F_k^{(i+1)}- F_k^{(i)}\right): 0 \le i < N_k\right\},
\end{equation}
and given $0 \le a < b  \le \infty$, what is the behavior of 
\begin{equation}\label{eq:gap:dist}
\lim_{k \rightarrow \infty} \frac{\left| G(k) \cap (a, b)\right|}{N_k}?
\end{equation}
If the sequence $F(k)$ is `truly random', that is, given by $$F(k) = \{ X_{(0)} \le X_{(1)} \le \ldots \le X_{(k)}\},$$ where the $\{X_{(i)}\}$ are the order statistics generated by independent, identically distributed (i.i.d.) uniform $[0, 1)$ random variables $\{X_n\}_{n=0}^{\infty}$, it is an exercise in probability theory to show that the gap distribution converges to a \emph{Poisson process} of intensity $1$. Precisely, for any $t>0$,
\begin{equation}\label{eq:pois}
\lim_{k \rightarrow \infty} \frac{\left| G(k) \cap (t, \infty)\right|}{N_k} = e^{-t}
\end{equation}
However, many sequences that arise `in nature' satisfy an equidistribution property but do not have Poissonian gaps. Following~\cite{Elkies}, we call such gap distributions \emph{exotic}. In this paper, we discuss in detail some examples of exotic gap distributions, which, moreover, can be calculated (or at least shown to exist) using methods arising from homogeneous dynamics, in particular dynamics of $SL(2, \R)$ actions on appropriate moduli spaces. In particular, the results we discuss share a similar philosophy; the sets $F(k)$ are associated to sets of angles or slopes of a discrete set of vectors in $\R^2$, and the gap distribution is studied by appropriate linear renormalizations, which can be viewed as part of an $SL(2, \R)$ action on an appropriate moduli space of geometric objects. The main novelty of this paper is the statement of three meta-theorems (Theorem~\ref{theorem:meta}, Theorem~\ref{theorem:meta:period}, and Theorem~\ref{theorem:meta:circle}), which give unified explanations of some of these examples by linking them to uniform distribution on various moduli spaces and which we expect can be used for future applications.

\subsection{Organization of the paper}\label{subsec:org} This paper is organized as follows: in  the remainder of this introduction we state results about our main (previously studied) examples: the Farey sequences $\F(Q)$ (\S\ref{subsubsec:farey}); slopes for lattice vectors (\S\ref{subsubsec:geomnum}); and saddle connection directions for translation surfaces (\S\ref{subsubsec:saddle}). We also briefly discuss the space of affine lattices and $\{ \{\sqrt n\} \}_{n \geq 1}$ in \S\ref{sec:affine}. In \S\ref{subsubsec:meta}, we state the main results Theorem~\ref{theorem:meta}, Theorem~\ref{theorem:meta:period}, and Theorem~\ref{theorem:meta:circle}. We describe how to use these results to explain our examples in \S\ref{sec:farey}-\S\ref{sec:saddle}, and prove the theorems in \S\ref{sec:metaproof}. Finally, in \S\ref{sec:quest}, we pose some natural questions suggested by our approach.

\subsection{Acknowledgements}\label{subsec:ack} Parts of this work are based on joint work with J.~Chaika, Y.~Cheung, and S.~Lelievre. I thank them for their stimulating intellectual partnership. Discussions with F.~Boca, M.~Boshernitzan, A.~Eskin, J.~Marklof, and A.~Zaharescu helped the author to see the common structural features in this circle of ideas, and without them this paper would not have been possible. In particular, the stunning paper~\cite{MS} of J.~Marklof and A.~Strombergsson has inspired not only this paper but the papers~\cite{AChaika, ACL, ACheung}. I would also like to thank the Polish Academy of Sciences for their hospitality, and in particular Piotr Przytycki for arranging a visit to Warsaw where many of the ideas in this paper were clarified. Finally, thanks are due to the organizers of the ICM Satellite Conference on Geometry, Topology, and Dynamics in Negative Curvature for the opportunity to speak at the meeting and also to contribute to this volume. We would also like to thank the anonymous referee for their careful reading and many remarks which clarified the exposition of this paper.

\subsection{Farey sequences}\label{subsubsec:farey} 
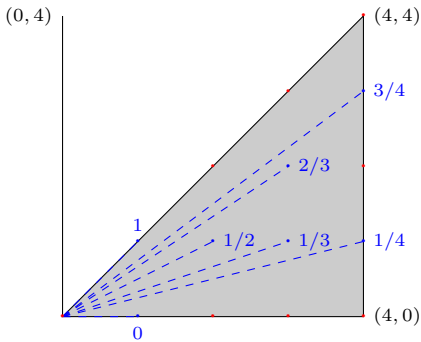
\begin{figure}\caption{The triangle $T_4$. Primitive integer vectors are given by dashed lines, and are labeled by their slopes.}\label{fig:triangle:image}
\begin{tikzpicture}[scale=1]
\filldraw[fill=black!20!white] (0,0)--(4,4)--(4,0)--cycle;
\draw (0,0)--(0,4);
\path (4, 0) node[right]{\tiny $(4,0)$};
\path (0, 4) node[left]{\tiny $(0,4)$};
\path (4, 4) node[right]{\tiny $(4,4)$};
\path (4, 4) node{\textcolor{red}{$\cdot$}};
\path (3, 3) node{\textcolor{red}{$\cdot$}};
\path (2, 2) node{\textcolor{red}{$\cdot$}};
\path (1, 1) node{\textcolor{blue}{$\cdot$}};
\path (4, 0) node{\textcolor{red}{$\cdot$}};
\path (3, 0) node{\textcolor{red}{$\cdot$}};
\path (2, 0) node{\textcolor{red}{$\cdot$}};
\path (1, 0) node{\textcolor{blue}{$\cdot$}};
\path (4, 3) node{\textcolor{blue}{$\cdot$}};
\path (4, 2) node{\textcolor{red}{$\cdot$}};
\path (4, 1) node{\textcolor{blue}{$\cdot$}};
\path (3, 2) node{\textcolor{blue}{$\cdot$}};
\path (3, 1) node{\textcolor{blue}{$\cdot$}};
\path (2, 1) node{\textcolor{blue}{$\cdot$}};
\draw[dashed, blue](0,0)--(1,1)node[blue, above]{\tiny $1$};
\draw[dashed, blue](0,0)--(4, 3)node[blue, right]{\tiny $3/4$};
\draw[dashed, blue](0,0)--(3,2)node[blue, right]{\tiny $2/3$};
\draw[dashed, blue](0,0)--(2,1)node[blue, right]{\tiny $1/2$};
\draw[dashed, blue](0,0)--(3,1)node[blue, right]{\tiny $1/3$};
\draw[dashed, blue](0,0)--(4,1)node[blue, right]{\tiny $1/4$};
\draw[dashed, blue](0,0)--(1,0)node[blue, below]{\tiny $0$};
\path (0, 0) node{\textcolor{red}{$\cdot$}};
\end{tikzpicture}
\end{figure}

Consider the integer lattice $\Z^2$. If we imagine an observer sitting at the origin $0$, the `visible' points in $\Z^2$ correspond to the set of \emph{primitive} vectors, that is, integer vectors which are not integer multiples of other integer vectors. If we consider slopes of vectors (as opposed to angles), it is natural to consider the set of vectors with slopes in $[0, 1]$. The set of slopes of (primitive) integer vectors with horizontal component at most $Q$ intersected with the interval $[0, 1]$  gives the \emph{Farey sequence} of level $Q$. More simply,  $\F(Q)$ consists of the set of fractions in between $0$ and $1$ with denominator at most $Q$. We write

$$\F(Q) : = \left\{ \gamma_0 = \frac 0 1 < \gamma_1 = \frac 1 Q < \gamma_2  \ldots < \gamma_i = \frac {p_i}{q_i} < \ldots \gamma_{N} = \frac 1 1\right\}$$
Here, $N = N(Q) = \sum_{i=1}^Q \varphi(i)$ is the cardinality of $\F(Q)$. By the above discussion, these correspond to the slopes of primitive integer vectors $\left(\begin{array}{c}q_i \\p_i\end{array}\right)$ in the (closed) triangle $T_Q$ with vertices at $(0,0)$, $(Q, 0)$, and $(Q, Q)$. That is, it is bounded above by the line $\{y=x\}$, below by the $x$-axis, and on the right by the line $\{x = Q\}$. The triangle $T_4$ is shown in Figure~\ref{fig:triangle:image}.

\begin{figure}[htbp]\caption{The region $A_{a,b}$ in inside the region $\Omega$.\medskip}\label{torusplot}
    \includegraphics[width=40mm, height = 40mm]{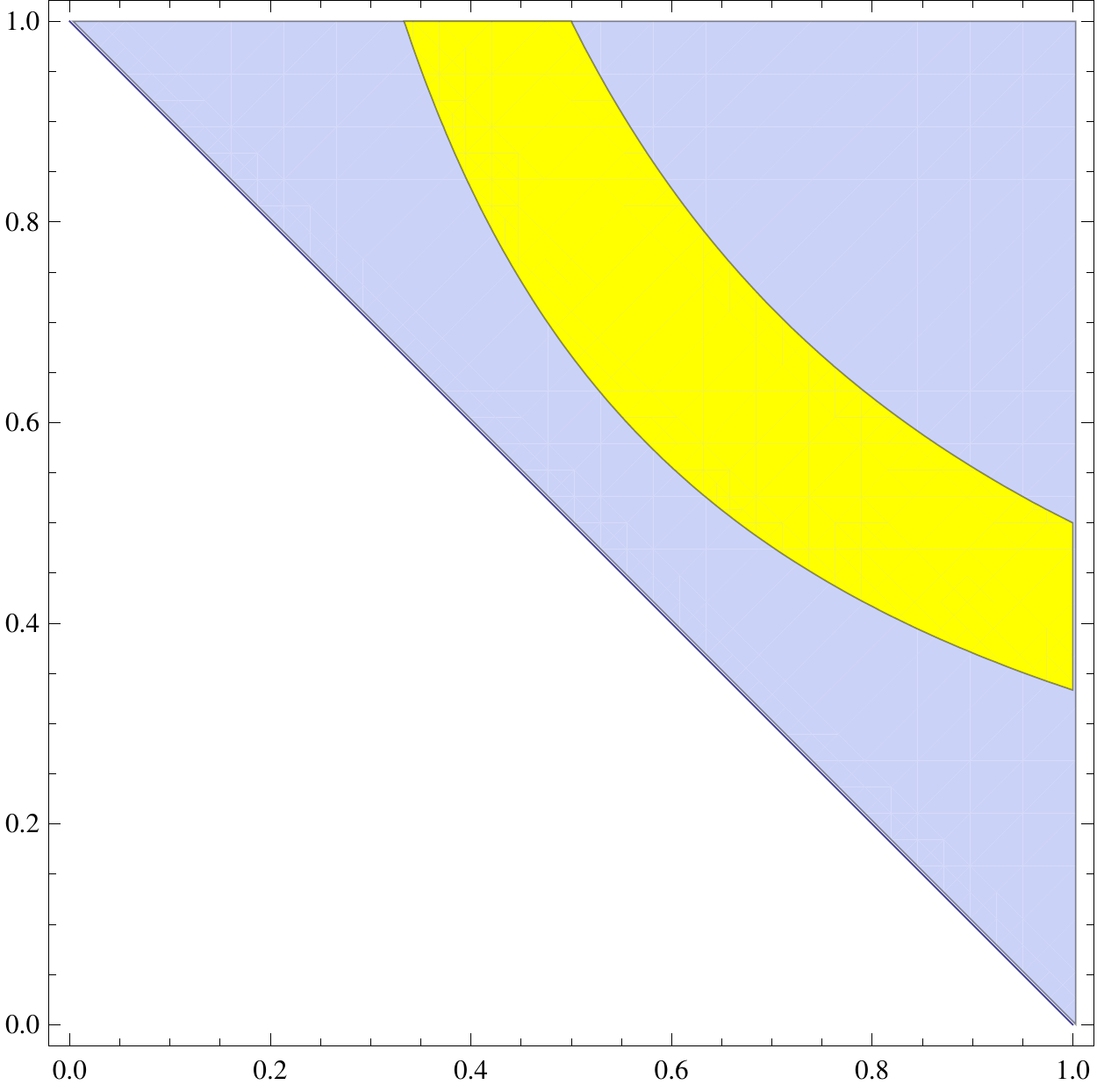}
    \end{figure}

\begin{figure}[htbp]\caption{The limiting distribution of gaps for Farey fractions and, (appropriately rescaled) lattice slopes.}\label{torusplot2}
  \includegraphics[width=0.6\textwidth, height = 40mm]{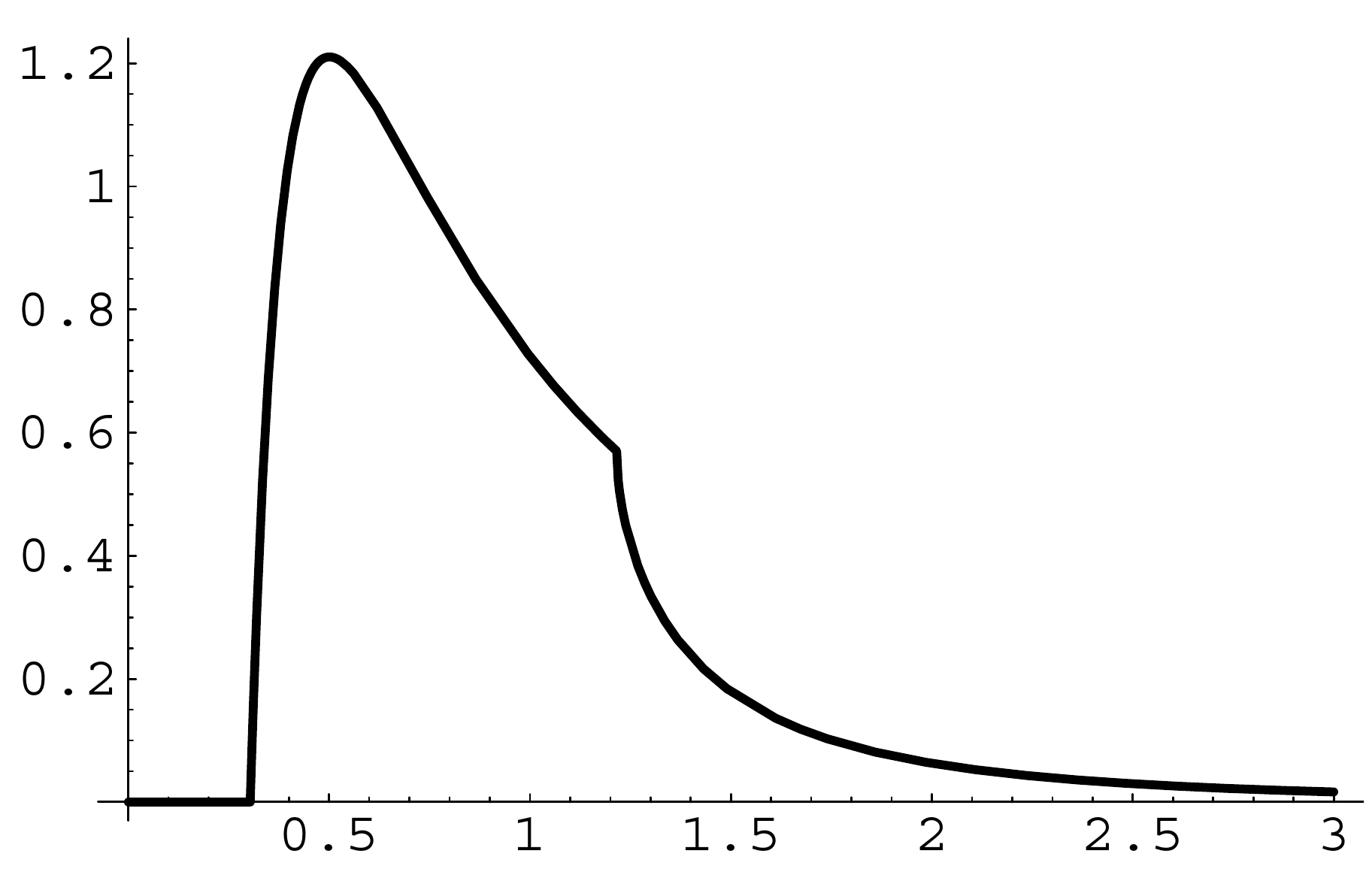}

\end{figure}

The sequences $\F(Q)$ equidistribute in $[0, 1]$ (by, for example, Weyl's criterion~\cite{Weyl}). We denote by $\mathcal{G}(Q)$ the set of normalized gaps between Farey fractions, that is, 
$$\mathcal{G}(Q) = \left\{N(Q)(\gamma_{i+1} - \gamma_i )= \frac{N(Q)}{q_i q_{i+1}}: 0 \le i < N(Q)\right \}.$$
The limiting distribution for $\mathcal{G}(Q)$ is given by the following beautiful theorem of R.~R. Hall, and illustrated in Figures~\ref{torusplot} and ~\ref{torusplot2}. Let $$\Om: = \{ (u, v) \in [0, 1]^2: u+v >1\},$$ and for $0 \le a, b < \infty$, let \begin{equation}\label{eq:aab}A_{a, b} = \left\{ (u,v) \in \Om:  b^{-1}< uv <  a^{-1}\right\},\end{equation} and set $\tilde{A}_{a,b} : = A_{\frac{\pi^2}{3}a, \frac{\pi^2}{3}b}$.
\medskip

\begin{Theorem}~\cite[R.R.Hall]{Hall}\label{theorem:hall} Fix $0\le a < b< \infty$. Then $$\lim_{Q \rightarrow \infty}\frac{ |\mathcal{G}(Q) \cap (a, b)|}{N(Q)}=2 |\tilde{A}_{a, b}|.$$\end{Theorem}

\noindent Differentiating the cumulative distribution function $$F_{Hall}(t) := |\tilde{A}_{0, t}|,$$ one can compute the probability distribution function $P_{Hall}(t)$ so that $$\int_{a, b} P_{Hall}(t) dt : = 2|\tilde{A}_{a, b}|.$$ We call this distribution (and any scalings) \emph{Hall's distribution}. The graph of $P_{Hall}(t)$ is given in Figure~\ref{torusplot2}, which is drawn from~\cite{BZsurvey}. The points of non-differentiability $\frac{3}{\pi^2}$ and $\frac{12}{ \pi^2}$ correspond to the transitions when the hyperbola $\left\{xy =    \frac{3}{\pi^2} t^{-1}\right\}$ enters the region $\Omega$  $\left(t=\frac{3}{\pi^2}\right)$ and when it hits the line $x+y=1$ $\left(t=\frac{12}{ \pi^2}\right)$.
In \S\ref{sec:farey}, we will, following~\cite{ACheung}, give a proof of Hall's theorem inspired by the work of F.~Boca, C.~Cobeli, and A.~Zaharescu~\cite{BCZ}. They created a map $T: \Omega \rightarrow \Omega$, now known as the \emph{BCZ map}, and used equidistribution properties of periodic orbits of this map to obtain many statistical results on $\F(Q)$ and $\mathcal{G}(Q)$. In~\cite{ACheung}, the author and Y.~Cheung showed that these results could be obtained by studying the horocycle flow on the space $X_2 = SL(2, \R)/SL(2, \Z)$ of unimodular lattices in $\R^2$. 

\subsubsection{Geometry of Numbers}\label{subsubsec:geomnum} One can also study the behavior of an arbitrary unimodular lattice $\La$.  Let $\La \subset \R^2$ be a unimodular lattice, and suppose $\La$ does not have vertical vectors. Let $\{s_1 < s_2 < \ldots < s_n< \ldots\}$ denote the slopes of the vectors (written in increasing order) in the vertical strip $V_1 =\{(u,v)^T: u \in (0, 1], v>0\}$. Here, and below, we use $(u, v)^T$ to denote the \emph{column vector} $\left(\begin{array}{c}u \\v\end{array}\right)$, as our matrices act on the left.  Let $$G_N(\La) = \{ s_{n+1} - s_n: 0 \le n \le N\}$$ denote the set of gaps in this sequence. Note that in this setting, we do not need to normalize, as our sequence is not contained in $[0, 1)$. Then we have that the limiting distribution of $G_N$ is also given by Hall's distribution. That is:

\begin{Theorem}\label{theorem:geomnum}\cite{ACheung} Let $0 \le a \le b \le \infty$. Then $$\lim_{N \rightarrow \infty} \frac{1}{N} | G_N(\La) \cap (a, b)| = 2|A_{a,b}|.$$
\end{Theorem}

\subsection{Saddle Connections}\label{subsubsec:saddle} We saw above that the Farey sequence could be interpreted geometrically as slopes of primitive integer vectors in $\R^2$. Primitive integer vectors also correspond to (parallel families) of closed geodesics on the torus $\R^2/\Z^2$, which can also be interpreted as closed billiard trajectories in the square $[0, 1/2]^2$. A natural generalization would be to try and understand similar families of trajectories for higher-genus surfaces, and/or for billiards in more complex polygons
More precisely, let $P$ be a Euclidean polygon with angles in $\pi\Q$. The billiard dynamical system on $P$ is given the (frictionless) motion of a point mass at unit speed with elastic collisions with the sides, satisfying the law of geometric optics: \emph{angle of incidence = angle of reflection}. A \emph{generalized diagonal} for the polygon $P$ is a trajectory for the billiard flow that starts at one vertex of $P$ and ends at another vertex. Since the group $\Delta_P$ generated by reflections in the sides of $P$ is finite, the \emph{angle} of a trajectory is well defined in $S^1 \cong S^1/\Delta_P$. The natural gap distribution question that arises in this context is: 

\begin{ques} What is the limiting distribution of the gaps between angles of generalized diagonals (normalized in terms of the length)?\end{ques}

More generally, one can ask about the limiting distribution for gaps for \emph{saddle connections}) in the more general setting of translation surfaces. A translation surface is a pair $(M, \om)$, where $M$ is a Riemann surface and $\om$ a holomorphic $1$-form. 

A saddle connection is a geodesic $\gamma$ in the flat metric induced by $\om$, connecting two zeros of $\om$. To each saddle connection $\gamma$ one can associate a holonomy vector $\vv_{\gamma} = \int_{\gamma} \om \in \C$. The set of holonomy vectors $\La_{sc}(\om)$ is a discrete subset of $\C \cong \R^2$, and varies equivariantly under the natural $SL(2, \R)$ action on the set of translation surfaces. Motivated by such concerns, and inspired by the work of Marklof-Strombergsson~\cite{MS} (of which more below in \S\ref{sec:affine}), the author and J.~Chaika~\cite{AChaika} studied the gap distribution for \emph{saddle connection directions}.  The relationship between flat surfaces and billiards in polygons is given by a natural \emph{unfolding} procedure, which associates to each (rational) polygon $P$ a translation surface $(X_P, \om_P)$. The main result of~\cite{AChaika} used the dynamics of the $SL(2,\R)$ action on the moduli space $\Omega_g$ of genus $g$ translation surfaces to show that generically, a limiting distribution exists. 

More precisely, given $R >0$, let 
\begin{equation}\label{eq:Theta:def} F^{\om}_R  : = \{\arg(\vv): \vv \in \Lambda_{\om} \cap B(0, R)\}\end{equation} denote the set of directions of saddle connections of length at most $R$. Masur~\cite{Masur} showed that the counting function $N(\omega, R): = |F^{\om}(R)|$ grows quadratically in $R$ for any $\omega$. Denote the associated normalized gap set by $G^{\om}(R)$.

\begin{Theorem}\label{theorem:sc:gap} (\cite[Theorem 1.1]{AChaika}) For almost every (with respect to Lebesgue measure on $\Omega_g$) translation surface $\om$, there is a limiting distribution for the gap set $G^{\om}(R)$. Moreover, this distribution has support at $0$, that is, for almost every $\omega \in \Omega_g$, and for any $\epsilon >0$, 
\begin{equation}\label{eq:main:proportion}
\lim_{R \rightarrow \infty} \frac{ |G^{\om}(R) \cap (0, \epsilon)|}{N(\omega, R)} >0.
\end{equation}
\end{Theorem}

\subsubsection{Lattice Surfaces}\label{subsubsec:veech}

The support at $0$ in Theorem~\ref{theorem:sc:gap} is in contrast to the setting of the torus, where, as seen in Figure~\ref{torusplot2}, there a gap between $0$ and $3/\pi^2$. This gap at $0$ is, in some sense, due to the symmetry of the torus- if we think of the $SL(2, \R)$ action on the moduli space $X_2$ of flat tori, the stabilizer of any point is (conjugate to) $SL(2, \Z)$. More generally,  It was shown in~\cite{AChaika} that if $\om$ is a \emph{lattice surface} (i.e., the stabilizer of the flat surface $\om$ under the $SL(2, \R)$ action is a lattice) that the limiting distribution for gaps has no support at $0$.

\begin{figure}\caption{The Golden $L$. The long sides of the $L$ each have length $\frac{1+ \sqrt{5}}{2}$. \medskip}\label{fig:golden:L}
\begin{tikzpicture}[scale=2]

\draw[red](0,0)--(1,0)node[black]{$\circ$};
\draw[blue] (1,0)--(1.6, 0)node[black]{$\circ$};
\draw[yellow](1.6, 0)--(1.6, 1)node[black]{$\circ$};
\draw[blue](1.6, 1)--(1,1)node[black]{$\circ$};
\draw[green](1,1)--(1, 1.6)node[black]{$\circ$};
\draw[red](1, 1.6)--(0, 1.6)node[black]{$\circ$};
\draw[green](0, 1.6)--(0, 1)node[black]{$\circ$};
\draw[yellow](0, 1)--(0,0)node[black]{$\circ$};

\draw[dashed](0,0)--(1.6, 1);
\draw[dashed](0,0)--(1, 1);
\draw[dashed](0,0)--(1, 1.6);
\draw[dashed](1,0)--(1.6, 1);

\path(1/2,0)node[below, red]{\tiny $1$};
\path(1.3,0)node[below, blue]{\tiny $\frac{\sqrt{5}-1}{2}$};
\path(0, 1/2)node[left, yellow]{\tiny $1$};
\path(0, 1.3)node[left, green]{\tiny $\frac{\sqrt{5}-1}{2}$};

\end{tikzpicture}

\end{figure}
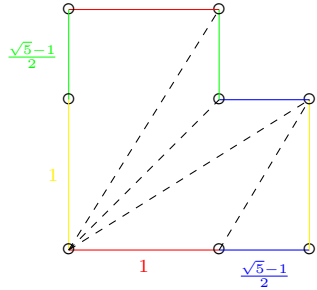

While it was in principle possible to compute the limiting distribution using the techniques in~\cite{AChaika}, the more geometric nature of the techniques in~\cite{ACheung} and the use of horocycle flows on moduli spaces can be generalized to the setting of lattice surfaces to give a roadmap for explicitly calculating the limiting distribution of gaps. In joint work~\cite{ACL} with J.~Chaika and S.~Lelievre, we proved Theorem~\ref{theorem:goldenL} on the gap distribution for the golden $L$, which is a surface of genus 2 with one double zero, displayed in Figure~\ref{fig:golden:L}. 

\begin{Theorem}\label{theorem:goldenL}\cite{ACL} There is an explicit limiting gap distribution for the set of slopes (equivalently, angles) for saddle connections on the golden $L$. The probability distribution function is differentiable except at a set of eight points. \end{Theorem}

\begin{figure}[htbp]\caption{The limiting and empirical distributions for gaps of saddle connection slopes on the golden $L$.\medskip}\label{fig:goldenLdist}
    \includegraphics[width=0.49\textwidth, height = 45mm]{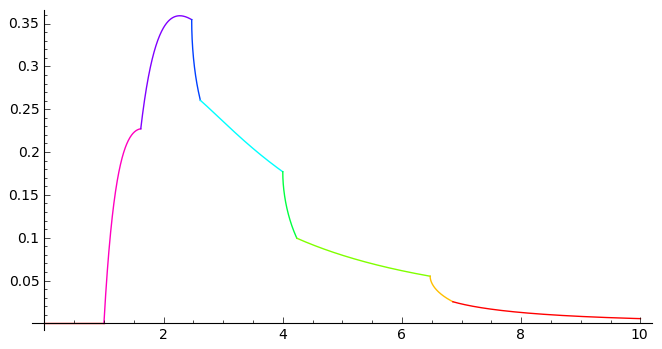}
  \includegraphics[width=0.49\textwidth, height = 45mm]{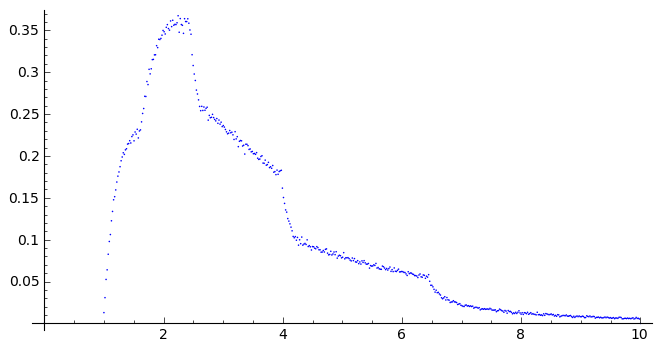}

\end{figure}

\noindent\textbf{Remark:} The limiting and empirical distributions are shown in Figure~\ref{fig:goldenLdist}, drawn from~\cite{ACL}. We refer the reader to ~\cite{ACL} for the precise formulas for the limiting distribution.

\subsection{Visible affine lattice points}\label{sec:affine} Another natural generalization of the Farey sequence is to consider \emph{affine lattices}, that is, translates of lattices by some fixed vector. We write $$\La = M \Z^2 + \vv,$$ where $M \in SL(2, \R)$ and $\vv \in \R^2$ (really $\vv$ is well-defined up to the lattice $M\Z^2$, so we think of it as an element of the torus $\R^2/M\Z^2$). Marklof-Strombergsson~\cite{MS} used dynamics on the space of affine lattices $\tilde{X}_2 = SL(2, \R)\ltimes \R^2/SL(2, \Z) \ltimes \Z^2$ to study the gap distribution for the angles of visible affine lattice points. They in fact considered much more general problems, studying the distribution of visible affine lattice points in higher dimensions, but for the purposes of this paper, we focus on their two-dimensional results.

\begin{figure}\caption{The wedge $A^{\theta}_R(\sigma).$\medskip }\label{fig:thinning:wedge}
\begin{tikzpicture}[scale=.8]
   \draw(0,0) circle (2.5cm);
  \draw (38:1)--(38:2.5);
  \path(0:2.6)node[right]{\tiny $R$};
      
  \draw(52:1)--(52:2.5);
  \filldraw[fill=black!20!white]  (0,0) -- (38:2.5) arc (38:52:2.5)--(52:3)node[above] {\tiny $(\theta +\frac{\sigma}{R^2})$} -- (0,0);
   \draw[dashed,red](0,0)--(45:3)node{\tiny\textcolor{red}{$\theta$}};
\draw(38:2.5)--(38:3)node[right] {\tiny $(\theta - \frac{\sigma}{R^2})$};

\end{tikzpicture}
\end{figure}
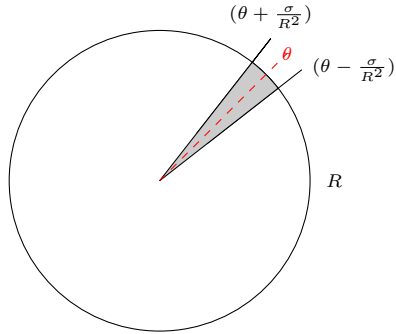

\noindent Consider the set of angles of lattice points of length at most $R$, that is, $$F_{\La} (R) : = \{ \arg \ww: \ww \in \La \cap B(0, R)\}.$$ To calculate the associated gap distribution $P_{\La}$, the key is to estimate the probability of finding multiple lattice points in `thinning' wedges.  Given $\sigma >0, \theta \in [0, 2\pi)$ and $R>0$ consider the wedge
$$A_R^{\theta} (\sigma) : =\{ \ww \in \R^2: \ww \in B(0, R), \arg(\ww) \in (\theta - \sigma R^{-2}, \theta + \sigma R^2)\},$$ shown in Figure~\ref{fig:thinning:wedge}. Here, the factor of $R^{-2}$ corresponds to the normalizing factor $\frac 1 N$ above, since the cardinality of $F_{\La}(R)$ is on the order of $R^2$. The gap distribution will be given by (the second derivative) of the limiting probability $$p_{\La, 0}(\sigma) = \lim_{R \rightarrow \infty} \lambda(\theta: A_R^{\theta} (\sigma) \cap \La = \emptyset)$$ that this wedge does not affine lattice points. This follows from the fact that if we let $P_{\La}(t)$ denote the probability distribution function of the limiting gap distribution, we have \begin{equation}\label{eq:gap:wedge} p_{\La, 0}(\sigma) = \sigma - \int_{0}^{\sigma} \int_{0}^{t} P_{\La}(s)ds dt\end{equation}

\noindent To compute $p_{\La, 0}(\sigma)$, note that rotating the region $A_R^{\theta} (\sigma)$ by the element \begin{equation}\label{eq:rtheta}r_{-\theta} = \left(\begin{array}{cc} \cos\theta& \sin\theta \\ -\sin\theta & \cos\theta
\end{array}\right),\end{equation} and scaling by \begin{equation}\label{eq:gt}g_t = \left(\begin{array}{cc} e^{-t/2} & 0 \\ 0 & e^{t/2}
\end{array}\right),\end{equation} (with $t = 2 \log R$) we obtain (approximately) the triangle $T(\sigma)$ with vertices at $(0, 0)$ and $(1, \pm \sigma)$. Thus, the probability that (for a randomly chosen $\theta$) the region $A_R^{\theta}$ does not contain points from $\La$ corresponds to the integral $$\int_{0}^{2\pi} h_{\sigma} (g_t r_{\theta} \La) d\theta ,$$ where $h_{\sigma}: \tilde{X_2} \mapsto \{0, 1\}$ is given by \begin{equation} h_{\sigma}(\La) =
\left\{ \begin{array}{ll}  1 & \La \cap T(\sigma) \neq \emptyset \\ 0 & \mbox{otherwise}\end{array}\right.\end{equation}
To understand the limits (as $R$, or equivalently $t \rightarrow \infty$) of the integrals $\int_{0}^{2\pi} h_{\sigma} (g_t r_{\theta} \La) d\theta$, we need to apply tools from Ratner's theory of unipotent flows. It can be shown that these integrals converge to $\int_{\tilde{X}_2} h_{\sigma} d\mu$, where the limit measure $\mu$ can be shown to be invariant under a unipotent subgroup $N \subset SL(2, \R) \ltimes \R^2$, following work of Shah~\cite{Shah}. Using Ratner's Theorem~\cite{Ratner}, the possibilities for the measure $\mu$ are either:
\medskip
\begin{itemize}
\item The Haar probability measure $\tilde{\mu}_2$
\medskip
\item A probability measure supported on the set of torsion points $\tilde{X_2}[n]$, that is, the support of the measure is restricted to affine lattices $M \Z^2 + \vv$ whose translation vector $\vv$ satisfies $n \vv \in M\Z^2$. 
\end{itemize}
\medskip

\noindent Depending on the properties of the initial affine lattice $\La$, torsion-supported measures can (and do) occur. As a consequence, the limiting gap distributions (which can all be explicitly computed) differ depending on the initial lattice, and in particular, whether the initial translation vector $\vv$ is a torsion point of the torus $\R^2/M \Z^2$ or not, and at what level~\cite[Theorem 1.3]{MS}. 

\subsubsection{$\sqrt{n}$ mod 1} Extraordinarily, the gap distribution for generic affine lattices $\La$ (that is, those for which the vector $\vv$ is not torsion, and so for which the limiting measure is Haar) coincides with the gap distribution for the sequence of fractional parts of $\{\sqrt{n}\}$. The reason they coincide is that both are coming from equidistribution of certain homogeneous trajectories with respect to the Haar measure on the space $\tilde{X_2}$. We strongly urge the reader to look at the beautiful paper~\cite{Elkies} for further details on the connection of $\sqrt{n}$ mod 1 to homogeneous dynamics. Perhaps ironically, despite the inspiration for our paper being provided from their pithy quote, our theorems below do not seem to give a direct explanation of their result.

\section{Meta-Theorems on Gap Distributions}\label{subsubsec:meta} As discussed in the introduction, the common thread that runs through these results is the creation of associated discrete sets of $\R^2$ and appropriate dynamical systems, which turns gap distribution questions into questions of equidistribution of orbits on certain moduli spaces. In this section, we state our main results, Theorem~\ref{theorem:meta}, Theorem~\ref{theorem:meta:period}, and Theorem~\ref{theorem:meta:circle}, in a very general setting, to attempt to explain several of these phenomena simultaneously. We first fix notation and our general framework in \S\ref{sec:setup}, before stating Theorem~\ref{theorem:meta} and Theorem~\ref{theorem:meta:period} in \S\ref{sec:theoremmeta:slope}, and Theorem~\ref{theorem:meta:circle} in \S\ref{sec:theoremmeta:angle}.

\subsection{Setup}\label{sec:setup} Our setup is very similar to that of~\cite{A}. Let $X$ be a locally compact metric space with a continuous $SL(2, \R)$-action, and a $SL(2, \R)$-equivariant assignment $$x \longmapsto \La_x$$ that associates to each point $x \in X$ a countable, discrete subset $\La_{x} \subset \R^2 \minuszero$. Equivariance, in this context, means that $$\La_{gx} = g\La_x$$ for all $g \in SL(2, \R), x \in X$, Our theorems will connect the gap distributions for slopes and angles in the set $\La_x$ to the orbit of the point $x$ under various subgroups of the $SL(2, \R)$ action.

\subsection{Slope gaps}\label{sec:theoremmeta:slope} The key idea of our meta-theorems comes from a simple observation on the behavior of slopes under (vertical) shears (a corresponding discussion can be made with inverse slopes and horizontal shears). Let $\uu = (u_1, u_2)^T, \vv = (v_1, v_2)^T \in \R^2$. Let $s_{\uu} = \frac{u_2}{u_1}$ and $s_{\vv} = \frac{v_2}{v_1}$ denote their slopes. Let \begin{equation}\label{eq:hs} h_s = \left(\begin{array}{cc}1 & 0 \\-s & 1\end{array}\right). \end{equation}

\begin{obs} For any $s \in \R$, he difference in slopes between $h_s \uu$ and $h_s \vv$ is the same as the difference in slopes between $\uu$ and $\vv$, that is, $$s_{h_s \uu} - s_{h_s \vv} = s_{\uu} - s_{\vv}.$$\end{obs}

\noindent This follows from the (even simpler) observation that for any $\uu \in \R^2$,  $$s_{h_s\uu} = s_{\uu} - s.$$ Thus, if we have a set of slopes $\s(\La) = \{s_1 < s_2 < \ldots s_n < \ldots\}$ of a family of vectors $\La \subset \R^2$, the behavior of the associated gap set $\G (\La) = \{ s_{i+1}-s_i: 1 \le i \le \infty\}$ is \emph{invariant} under the action of $h_s$, that is $$\G(h_s \La) = \G(\La).$$ In particular, by considering the flow $h_s$ at the times $s_1, s_2, \ldots$, we can recover the gap $s_{i+1}-s_i$ by sampling the smallest positive element of $\s(h_{s_i} \La)$. 

\subsubsection{Minkowski properties, short vectors, and vertical trips}\label{subsec:mink} For our theorem on slopes, we need two additional technical assumptions on the sets $\La_x$, a few definitions.

\medskip
\noindent\textbf{Symmetry.} $\La_x$ is centrally symmetric, that is $-\La_x = \La_x$.
\medskip

\noindent\textbf{Minkowski Property.} For all $x \in X$, there is a $c = c(x)$ such that for any convex, centrally symmetric set $K \subset \R^2$ of volume at least $c$, $$K \cap \La_x \neq \emptyset.$$
We leave it as an exercise to the reader that the minimal such constant $c(x)$ is \emph{invariant} under the action of $SL(2, \R)$.
\medskip

\noindent\textbf{Vertically and horizontally short points.} We say that $x \in X$ is \emph{vertical} (respectively \emph{horizontal}) if $\La_x$ contains a vertical (resp. horizontal) vector. Given a constant $\eta$, we say that $x$ is $\eta$-\emph{vertically short} (resp. $\eta$-\emph{horizontally short}) if $\La_x$ contains a vertical (resp. horizontal) vector of length at most $\eta$. 
\medskip

\noindent\textbf{Vertical Strip Property.}  For $\eta>0$, et $V_{\eta} : = \{(u,v)^T: 0 < u \le \eta, v > 0\}$ denote the vertical strip of width $\eta$ in $\R^2$. Suppose $x$ is horizontal. Then $\La_x \cap V_{\eta}$ is non-empty for any $\eta >0$.

\medskip

\noindent\textbf{Exceptional points} We say that $x$ is $\eta$-\emph{exceptional} if it is $\frac{\eta}{4c(x)}$-vertically short, that is, if it contains a vertical vector of length at most $\frac{\eta}{4c(x)}$. We denote by $X_{\eta}$ the set of non-$\eta$-exceptional points.

\subsubsection{Horocycles and transversals} Note that since $h_s$ fixes vertical vectors, $X_{\eta}$ is an $h_s$-invariant set for any $\eta >0$. An important lemma, crucial to the statement of Theorem~\ref{theorem:meta} is:

\begin{lemma}\label{lemma:omt} Let $\eta>0$, and let $\Om_{\eta} \subset X$ denote the set of $\eta$-horizontally short points. Then $\Om_{\eta}$ is a transversal to the $h_s$ action on the set $X_{\eta}$. That is, for every $x \in X_{\eta}$, the $h_s$-orbit of $x$ intersects $\Om_{\eta}$ in a non-empty, discrete set of times. \end{lemma}
\medskip
\noindent We prove this lemma in \S\ref{sec:metaproof}. The crucial consequence of this lemma, for us, is the existence of the induced return map $T_{\eta}: \Om_{\eta} \rightarrow \Om_{\eta}$ and return time function $R_{\eta}: \Om_{\eta} \rightarrow \R^+$, given by $$R_{\eta}(x) = \min\{s>0: h_s x \in \Om_{\eta}\}$$ and $T^t(x) = h_{R(x)} x.$ In particular, any finite $h_s$-invariant measure $\nu$ on $X_{\eta}$ can be decomposed as $d\nu = d\tilde{\nu} ds$, where $\tilde{\nu}$ is a $T_{\eta}$-invariant probability measure on $\Om_{\eta}$ and we identify $X_{\eta}$ with the suspension space $$ \{ (x, s): x \in \Om_{\eta}, 0 \le s \le R_{\eta}(x)\}/\sim,$$ where $(x, R_{\eta}(x)) \sim (T_{\eta}(x), 0)$.

\subsubsection{Long horocycles}
Define \begin{equation}\label{eq:slopeset} \s_{\eta}(x) : = \{ 0 \le s_1^{\eta}(x)< s_2^{\eta}(x) < \ldots < s_n^{\eta}(x) < \ldots \} \end{equation} Let \begin{equation}\label{eq:gapslope} \G_N^{\eta}(x):=\{1 \le i \le N: s_{i+1}^{\eta}(x) - s_i^{\eta}(x)\}\end{equation} denote the set of the first $N$ slope gaps. Here, we do not normalize, since the slopes $s_i^t$ are going to $\infty$ in all of our applications. Given $x \in X$, and $S>0$, let $\sigma_{x, S}$ denote the Lebesgue probability measure $d\sigma_{x, S} = \frac{1}{S}ds$ on the segment $\{h_s x: 0 \le s \le S\}.$

\begin{Theorem}\label{theorem:meta} Let $t>0$. Suppose $\sigma_{x,S} \rightarrow \nu$ as $S \rightarrow \infty$ (here, and below, convergence of measures is in the weak-$\star$ topology), where $\nu$ is a finite measure supported on $X_{\eta}$, and $d\nu = d\tilde{\nu} ds$ as above. Then for any $0 \le a \le b \le \infty$, we have \begin{equation}\label{eq:slopelim} \lim_{N \rightarrow \infty} \frac{| \G_N^{\eta}(x) \cap [a, b]|}{N} = \tilde{\nu}\left( R_{\eta}^{-1}([a,b])\right).\end{equation} That is, the slope gap distribution is given by the distribution of the return time function $R_{\eta}$.

\end{Theorem}
\noindent
\medskip This theorem has the following immediate corollary:

\begin{Cor}\label{cor:meta:ergodic} Let $\nu$ be a ergodic, $h_s$-invariant probability measure supported on $X_{\eta}$. Then for $\nu$-almost every $x$, $$\lim_{N \rightarrow \infty} \frac{| \G_N^{\eta}(x) \cap [a, b]|}{N} = \tilde{\nu}\left( R_{\eta}^{-1}([a,b])\right).$$
\end{Cor}

\subsubsection{Long closed horocycles} We can also state a theorem for $x$ with periodic $h_s$-orbits. Again, we need a definition and a preliminary (elementary) lemma. We say that $x \in X$ is $h_s$-\emph{periodic} with period $s_0>0$ if $h_{s_0} x =x$, and $h_s x \neq x$ for all $0 < s < s_0$.

\begin{lemma}\label{lemma:period} Let $x$ be $h_s$-periodic. Then for every $\eta>0$ there is an $N_0 = N_0(\eta)$ such that for any $N \geq N_0$, $$\G_{N}^{\eta}(x) = \G_{N_0}^{\eta}(x).$$ That is, the set of gaps repeats.\end{lemma}
\noindent We postpone the proof of the lemma to \S\ref{sec:metaproof}. Let $g_t= \mbox{diag}(e^{t/2}, e^{-t/2})$ be as in (\ref{eq:gt}). Then we have the conjugation relation \begin{equation}\label{eq:conj} g_t h_s g_{-t} = h_{se^{-t}}. \end{equation} 
Thus, if $x$ is $h_s$-periodic with period $s_0$, $g_{-t}x$ is $h_{s_0 e^{t}}$-periodic, since $$h_{s_0e^t}g_{-t} x = g_{-t} h_{s_0} x = g_{-t} x.$$ Note also that if $x \in X_{\eta}$, then $g_{-t} x \in X_{e^{t/2} \eta}$ for any $t>0$, since $g_{-t}$ expands by $e^{t/2}$ in the vertical direction. Let $\rho_{x, T}$ denote the Lebesgue probability measure $d\rho_{x, T} = \frac{ds}{e^T s_0}$ supported on the periodic orbit $$\{h_s g_{-T} x: 0 \le s \le e^T s_0\}.$$

\begin{Theorem}\label{theorem:meta:period} Suppose $\rho_{x, T} \rightarrow \nu$ as $T \rightarrow \infty$, with $\nu$ a finite measure supported on $X_{\eta_0}$ for some $\eta_0 >0$ (in fact, it will be supported on any $X_{\eta}$), and $d\nu =  d \tilde{\nu} ds$, with $\tilde{\nu}$ the $T_{\eta_0}$-invariant probability measure on $X_{\eta_0}$. Then \begin{equation}\label{eq:slopeperiod} \lim_{\eta \rightarrow \infty} \frac{\left| \left(\frac{\eta}{\eta_0}\right)^2 \G_{N_0(\eta)}^{\eta}(x) \cap [a, b]\right|}{N_0(\eta)} = \tilde{\nu}\left( R_{\eta_0}^{-1}([a,b])\right).\end{equation} 
To re-iterate, the right hand choice of this limit does not depend on the choice of $\eta_0$, which can be chosen to be any positive number.\end{Theorem}

\subsection{Angle gaps and thinning wedges}\label{sec:theoremmeta:angle} In this section, we work with angles, as opposed to slopes. We also do not need the assumptions of central symmetry and the Minkowski property, though they are present in our applications. Let $x \in X$, and let $$\Theta_x(R) : = \{0 < \theta_1 < \theta_2 < \ldots < \theta_N < 2\pi\}$$ denote the set of angles of vectors in the set $\La_x \cap B(0, R)$, where $B(0, R) =\{\vv \in \R^2: \|\vv\|_2 \le R\}$ is the ball of radius $R$ in $\R^2$, and $N = N(R)$ is the cardinality of the set of distinct angles. Arguing as in \S\ref{sec:affine}, the gap distribution for $\Theta_x(R)$ can be reduced to studying the limiting probability $$p_{0}(x, \sigma) = \lim_{R \rightarrow \infty} \lambda(\theta: A_{\theta}^{\sigma}(R) \cap \La_x = \emptyset),$$ where $\lambda$ denotes the Lebesgue probability measure on $S^1$, $\sigma>0$ is a fixed parameter, and \begin{equation} A_R^{\theta} (\sigma) : =\{ \vv \in \R^2: \vv \in B(0, R), \arg(\vv) \in (\theta - \sigma R^{-2}, \theta + \sigma R^2)\},\end{equation} as shown in Figure~\ref{fig:thinning:wedge}. Let $\lambda_{x, R}$ denote the Lebesgue probability measure $d\lambda_{x, R} = \frac{d\theta}{2\pi}$ supported on the orbit $\{g_t r_{\theta} x \}_{0 \le \theta < 2\pi},$ with $t = 2\ln R$.

\begin{Theorem}\label{theorem:meta:circle} Suppose $\lambda_{x, R} \rightarrow \nu$ as $R \rightarrow \infty$, where $\nu$ is a probability measure on $X$. Then \begin{equation}\label{eq:circle:limit} p_0(x, \sigma) = \nu(x \in X: \La_x \cap T(\sigma) = \emptyset), \end{equation} where $T(\sigma)$ is the triangle with vertices at $(0, 0)$ and $(1, \pm \sigma)$. Moreover, if we define $$p_{i}(x, \sigma) : = \lim_{R \rightarrow \infty} \lambda(\theta: |A_{\theta}^{\sigma}(R) \cap \La_x| = i),$$ we have \begin{equation}\label{eq:circle:limit:i}p_i(x, \sigma) = \nu(x \in X: |\La_x \cap T(\sigma)| = i).\end{equation}

\end{Theorem}

\section{Farey sequences}\label{sec:farey} \noindent In this section, we show how to use Theorem~\ref{theorem:meta:period} to prove Hall's Theorem~\ref{theorem:hall}, and how to use Theorem~\ref{theorem:meta} to prove the geometry of numbers result Theorem~\ref{theorem:geomnum}. We follow the exposition in~\cite{ACheung}, where it was shown how these theorems can be obtained from results on equidistribution of long orbits of the group $\{h_s\}$ on the space of unimodular lattices $X_2 = G/\Gamma$, where $G = SL(2, \R)$ and $\Gamma = SL(2, \Z)$. Following the notation in \S\ref{subsubsec:meta}, our space $X$ is given by the space $X_2$, the assignment of a discrete set is given by $$gSL(2, \Z) \longmapsto g \Z^2_{prim},$$ and we see immediately that the Minkowski and symmetry conditions of \S\ref{subsec:mink} are satisfied. It is also possible to verify the vertical strip condition explicitly (see~\cite{ACheung}). We note that $\Z^2$ is $h_s$-periodic with period $1$, since $h_{1} \Z^2 = \Z^2$.  Given $\eta>0$, the set of slopes in $V_{\eta}$ is the set of non-negative rational numbers with denominator at most $Q = \lfloor \eta \rfloor$. That is, we can write $$\s_{N}^{\eta}(\Z^2) : = \{n+\F(Q): n \geq 0 \},$$ where $n$ ranges over the non-negative integers. Thus, the gap set $\G_N^{\eta}(\Z^2)$ corresponds to the set of gaps in the Farey sequence $\F(Q)$, and has cardinality $N(Q)$. We recall that Sarnak~\cite{Sarnak} proved that the measures $\rho_{\Z^2, T}$ equidistribute with respect to the Haar measure probability $\mu_2$. That is, $$\rho_{\Z^2, T} \rightarrow \mu_2.$$ Thus, to complete the proof of Theorem~\ref{theorem:hall} (assuming Theorem~\ref{theorem:meta:period}), we need to describe the distribution of the return time function (with respect to the disintegrated measure) on the transversal $$\Om_{\eta} :=\{\La \in X_2: \La \mbox{ is } \eta-\mbox{horizontally short}\},$$ for some choice of $\eta>0$. 

\subsection{Description of transversal and return map}\label{subsec:transversal} The main result of~\cite{ACheung} is to explicitly describe the transversal $\Om_1$ (more generally, $\Om_{\eta}$) and the associated return map and return time function. It was shown that $\Om_{\eta}$ is in bijective correspondence with the set $$\{(a, b) \in \R^2: a, b \in (0, \eta], a+b > \eta\},$$ via the map $$(a, b) \mapsto \La_{a,b} = \left(\begin{array}{cc}a & b \\0 & a^{-1}\end{array}\right)\Z^2$$ and that the roof function is given by $$R_{\eta}(a, b) = \frac{1}{ab}.$$ In particular, it does not depend on $\eta$. If we set $\eta = 1$, the return map is the \emph{BCZ map} (named for its creators, Boca-Cobeli-Zaharescu~\cite{BCZ}), given by $$T_1(a, b) = \left(b, -a + \left\lfloor \frac{1+a}{b} \right\rfloor\right).$$ More generally, $$T_{\eta}(a, b) = \left(b, -a + \left\lfloor \frac{\eta+a}{b} \right\rfloor\right).$$  $T_{\eta}$ is piecewise linear, and all the linear maps have determinant 1, and so the Lebesgue probability measure $\frac{2}{\eta^2} dadb$  on $\Omega_b$ is $T_{\eta}$-invariant, and that $dadbds$ is, up to scaling, Haar measure on $X_2$ (which gives full measure to the set $X_{2, \eta}$ for any $\eta >0$). The distribution of the roof function $R$ can be seen to be independent of $\eta$ (see \S\ref{sec:selfsim} for an explanation of a general self-similarity phenomenon), so we can choose $\eta =1$, obtaining the result of Theorem~\ref{theorem:hall} up to the natural normalization factor $\frac{3}{\pi^2}$, which occurs since in the statement of Theorem~\ref{theorem:hall} we normalize the gaps with a factor of $N(Q) \sim \frac{3}{\pi^2} Q^2$, as opposed to the factor $Q^2$ in the statement of Theorem~\ref{theorem:meta:period}.\qed

\subsection{Proof of Theorem~\ref{theorem:geomnum}} The above computations are also crucial in the proof of Theorem~\ref{theorem:geomnum} assuming Theorem~\ref{theorem:meta}. The other crucial step in the proof is Dani's measure classification~\cite{Dani} result for $h_s$-invariant measures on $X_2$, which states that any ergodic $h_s$-invariant probability measure on $X_2$ is either supported on a periodic orbit or Haar measure. A lattice is $h_s$-periodic if and only if it has a vertical vector. To show that the existence of a vertical vector implies periodicity, note that such lattices can be written as $$\La = \left(\begin{array}{cc}0 & -a^{-1} \\a & b\end{array}\right)\Z^2,$$ and it is a direct verification to show $$h_{a^2}\La = \left(\begin{array}{cc}0 & -a^{-1} \\a & b-a\end{array}\right)\Z^2 = \left(\begin{array}{cc}0 & -a^{-1} \\a & b-a\end{array}\right)\left(\begin{array}{cc}1 & 1 \\0 & 1\end{array}\right)\Z^2 = \La.$$ On the other hand, if $\La$ is $h_s$ periodic, it must be \emph{divergent} under $g_t$, which is equivalent to the existence of a vertical vector. This measure classification and characterization, combined with the Birkhoff ergodic theorem, yields:

\begin{lemma}\label{lem:equidist:vert} Let $\La \in X_2$ be a lattice without vertical vectors. Then the measures $\sigma_{\La, S}$ converge to the Haar probability measure $\mu_2$. \end{lemma}
\medskip
\noindent Thus, applying Theorem~\ref{theorem:meta}, we see that for any $\eta >0$, $$\lim_{N \rightarrow \infty} \frac{1}{N} | \G^{\eta}_N(\La) \cap (a, b)| = \frac{2}{\eta^2}|\{(u,v) \in \Om_{\eta}: R_{\eta}(u,v) \in (a, b)\}| = 2|A_{a,b}| $$ where $A_{a,b}$ is as in (\ref{eq:aab}), and $| \cdot |$ is Lebesgue measure $dxdy$. The statement of Theorem~\ref{theorem:geomnum} is the above statement for $\eta =1$. \qed

\section{Saddle connections}\label{sec:saddle} In this section, we sketch a proof of Theorem~\ref{theorem:sc:gap}, following a similar strategy to the outline in \S\ref{sec:affine}, and attempting to indicate how it can be seen as a consequence of Theorem~\ref{theorem:meta:circle}. We will particularly focus on the conclusion that for generic tranlsation surfaces, the gap distribution has support at $0$. Recall notation: $\Om_g$ is the moduli space of holomorphic $1$-forms on compact genus $g$ Riemann surfaces, and a saddle connection $\gamma$ on a surface $(M, \om) \in \Om_g$ is a geodesic (with respect to the flat metric given by the one-form) connecting two zeros of $\om$ with no zeros in its interior. The holonomy vector $\vv_{\gamma} \subset \C \equiv \R^2$ associated to $\gamma$ is given by $$\vv_{\gamma} = \int_{\gamma} \om.$$ 
The set of holonomy vectors $\La_{\om}$ is a discrete subset of $\R^2$, and we are interested in the set of \emph{small gaps} between directions of vectors. Precisely, letting $F^{\om}_R  : = \{\arg(\vv): \vv \in \Lambda_{\om} \cap B(0, R)\}$ denote the set of directions, we are interested in the normalized gap set $G^{\om}(R)$, and in particular the size of the set $G^{\om}(R) \cap (0, \epsilon)$ for $\epsilon >0$. Note that having a small gap is equivalent to having two directions in a `thin wedge'.

\subsection{$SL(2, \R)$-action}\label{sec:sl2r} A point $(M, \om) \in \Om_g$ determines (and is determined by) an atlas of charts from the surface punctured at the zeros of $\om$ to the plane $\C$, given by integration of $\om$. The transition maps for these charts are translations, and $\om$ is given by $dz$ in these coordinates. This also gives a natural $SL(2, \R)$ action on $\Om_g$ via linear postcomposition with charts. There is (\cite{Masur:IET},~\cite{Veech}) a natural $SL(2, \R)$-invariant probability measure on $\Om_g$, which we denote $\mu_g$. 

\subsubsection{Equivariance} The assigment $(M, \om) \mapsto \La_{\om}$ is an $SL(2, \R)$-equivariant assignment, in the sense that for any $h \in SL(2, \R)$, $$\La_{ h \om} = h \La_{\om}.$$

\subsection{Renormalization} As above, we consider the wedges $A_R^{\theta}(\sigma)$. To understand the support of the distribution at $0$, we consider probability that for a randomly chosen $\theta$ there are at least $2$ points in the thinning wedge $A_R^{\theta}(\sigma)$. This will give us the probability that our gap set has gaps of size at most $\sigma$. As above, by renormalizing by $g_t$ and $r_{\theta}$, we can write this probability as 
$$\int_{0}^{2\pi} h_{2, \sigma} (g_t r_{\theta} \La) d\theta ,$$ where $h_{2, \sigma}: \Om_g \mapsto \{0, 1\}$ is given by \begin{equation} h_{2, \sigma}(\om) =
\left\{ \begin{array}{ll}  1 & |\La_{\om} \cap T(\sigma)| \geq 2 \\ 0 & \mbox{otherwise}\end{array}\right.\end{equation}

\subsection{Equidistribution} Since for $g \geq 1$ the space $\Om_g$ is not a homogeneous space, one cannot use Ratner-type technology to classify the possible limit measures for the measures $\lambda_{\om, R}$ supported on $\{g_t r_{\theta} \om\}_{0 \le \theta <2\pi}$. However, Eskin-Masur~\cite[Theorem 1.5]{EM} used a general ergodic result of Nevo to show that for $\mu_g$-almost every starting point $\om$, the measures converge (in an appropriately smoothed sense) to $\mu_g$, that is, the main assumption of Theorem~\ref{theorem:meta:circle} is satisfied, with $\nu = \mu_g$. It is worth noting that \emph{there is no known example of a point $\om \in \Om_g$ for which there does not exist a measure $\nu$ satisfying the hypothesis of Theorem~\ref{theorem:meta:circle}. } Using Nevo's theorem, the author and J.~Chaika~\cite{AChaika} were able to conclude:

$$\lim_{t \rightarrow \infty} \int_{0}^{2\pi} h_{2, \sigma} (g_t r_{\theta} \La) d\theta = \int_{\Om_g} h_{2, \sigma} d\mu_g,$$
Finally, a measure estimate on the set of surfaces with simultaneous short saddle connections due to~\cite{MasurSmillie} (see also~\cite[Lemma 7.1]{EMZ}) suffices to show that the integral on the right hand side is positive for any $\sigma >0$. 

\subsection{Slopes} For certain measures, Theorem~\ref{theorem:meta} and Corollary~\ref{cor:meta:ergodic} be applied in this setting, as the Minkowski condition is satisfied by~\cite[Theorem 1]{HS}. However, to give gap distributions for slopes of surfaces in the support of $h_s$-invariant measures on $\Om_g$, one still has to verify the vertical strip condition for the support of these measures, and then, identify the transversal, the roof function, and the associated measure $\tilde{\nu}$ on the transversal in order to get explicit formulas. If we have a \emph{closed} $SL(2, \R)$ orbit, then the points in that orbit are lattice surfaces, where it is possible to carry out this program.

\begin{figure}[htbp]\caption{The transversal $X_1$ to the horocycle flow $h_s$ for the $SL(2, \R)$-orbit of the golden $L$. The colors indicate different behaviors of the return map.\medskip}\label{fig:goldenLtrans}
\includegraphics[width=32mm, height =50mm]{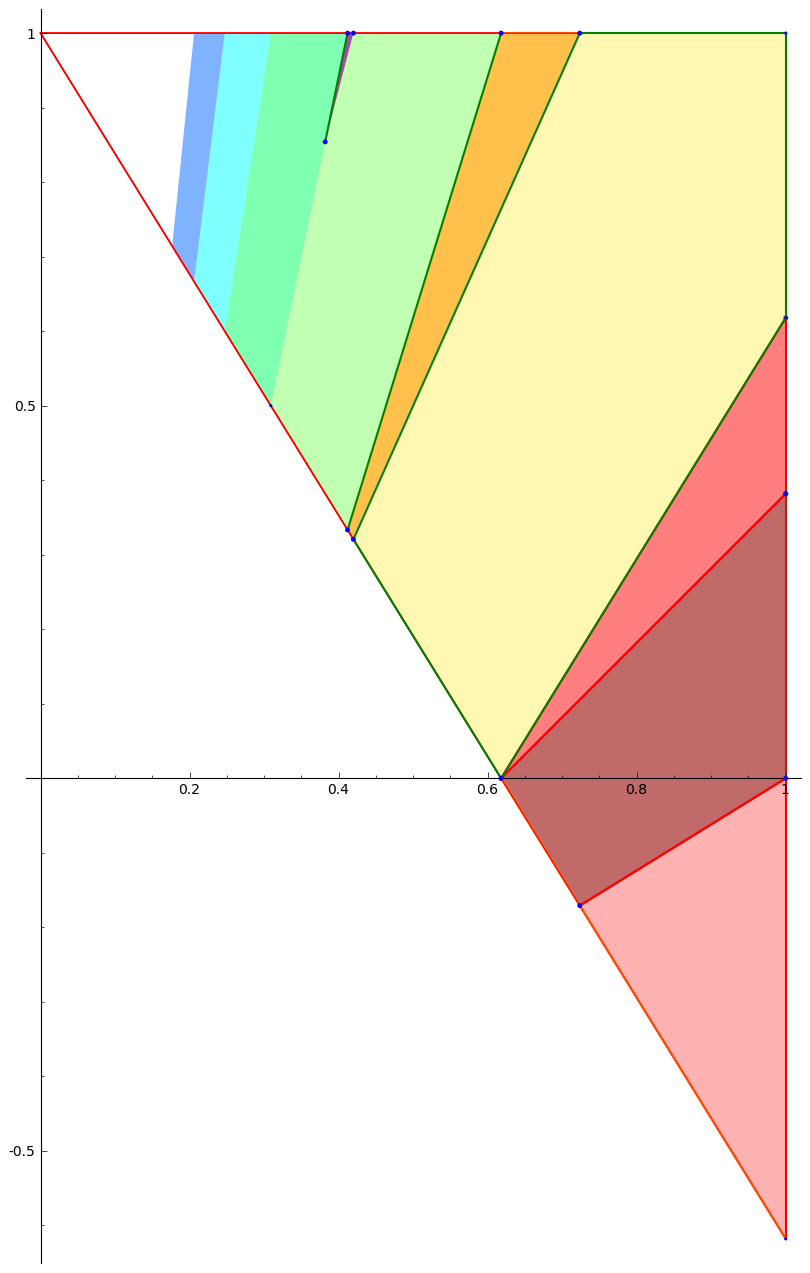}
\end{figure}

\subsection{Lattice Surfaces} If $(M, \om)$ is a lattice surface (for example the golden $L$), with stabilizer $\Gamma \subset SL(2, \R)$, we restrict our parameter space $X$ to from $\Om_g$ to the $SL(2, \R)$-orbit of $(M, \om)$, which can be identified with $SL(2, \R)/\Gamma$ via $$g(M, \om) \longmapsto g\Gamma.$$ It is possible (see~\cite{ACL}) to verify the vertical strip condition in this setting. Since the stabilizer $\Gamma$ must be non-uniform, we can assume that $(M, \om)$ is $h_s$-periodic, that is, that $\Gamma$ contains an element $h_{s_0}$ for some $s_0 >0$. Applying Theorem~\ref{theorem:meta:period}, we obtain that the limiting distribution of the slope gaps for vectors in $\La_{\om}$ (in the first quadrant) is given by the distribution of the return time function $R$ for $h_s$ to the transversal $X_1 : = \{ g\Gamma \in X: \La_{g \om} \mbox{ is $1$-horizontally short } \}.$ In~\cite{ACL}, the transversal $X_1$, the return time function, and the return map were explicitly computed for the translation surface $(M, \om)$ associated to the golden $L$. The transversal, in this case, can be identified with the triangle with vertices at $(0, 1)$, $(1, 1)$, and $\left(1, \frac{1-\sqrt 5}{2}\right)$, shown in Figure~\ref{fig:goldenLtrans}. In principal, with some effort, this program can be carried out for any lattice surface. 
%

\section{Abstract Homogeneous Dynamics}\label{sec:metaproof} \noindent In this section, we prove our main results Theorem~\ref{theorem:meta}, Theorem~\ref{theorem:meta:period}, and Theorem~\ref{theorem:meta:circle}. We prove Theorem~\ref{theorem:meta} and Theorem~\ref{theorem:meta:period} in \S\ref{subsec:horo:meta} and Theorem~\ref{theorem:meta:circle} in \S\ref{subsec:circle:meta}.

\subsection{Horocycle Flows}\label{subsec:horo:meta} We recall notation: $X$ is our parameter space, equipped with an $SL(2, \R)$-action, and given $x \in X$ we have the associated discrete set $\La_x \subset \R^2$, satisfying:

\medskip

\begin{description}

\item[Equivariance] For any $g \in SL(2, \R)$, $\La_{gx} = g\La_x$.

\medskip
\item[Symmetry] $\La_x$ is centrally symmetric, that is $-\La_x = \La_x$

\medskip

\item[Minkowski Property] For all $x \in X$, there is a $c = c(x)$ such that for any convex, centrally symmetric set $K \subset \R^2$ of volume at least $c$, $$K \cap \La_x \neq \emptyset.$$
\medskip

\item[Vertical Strip Property.]  Fix $\eta>0$. Let $V_{\eta} : = \{(u,v)^T: 0 < u \le \eta, v > 0\}$ denote the vertical strip in $\R^2$. Suppose $x$ is horizontal. Then $\La_x \cap V_{\eta}$ is non-empty.\medskip
\end{description}

\noindent Recall that we say that $x \in X$ is \emph{vertical} (respectively \emph{horizontal}) if $\La_x$ contains a vertical (resp. horizontal) vector, and that given a constant $\eta$, we say that $x$ is $\eta$-\emph{vertically short} (resp. $\eta$-\emph{horizontally short}) if $\La_x$ contains a vertical (resp. horizontal) vector of length at most $\eta$. We say that $x$ is $\eta$-\emph{exceptional} if it is $\frac{\eta}{4c(x)}$-vertically short, that is, if it contains a vertical vector of length at most $\frac{\eta}{4c(x)}$. We denote by $X_{\eta}$ the set of non-$\eta$-exceptional points. Note that $X_{\eta}$ and $X \backslash X_{\eta}$ are both $h_s$-invariant sets.

\subsubsection{Existence of return map}\label{subsubsec:return} In this subsection, we prove Lemma~\ref{lemma:omt}, which gives the existence of a well-defined transversal $\Om_{\eta}$ and return map $R_{\eta}$ for the flow $h_s$ on the invariant set $X_{\eta}$, for any fixed $\eta >0$. Recall that $$\Om_{\eta} = \{x \in X: x \mbox{ is } \eta-\mbox{horizontally short}\},$$ and $$V_{\eta} =\{ (u, v)^T \in \R^2: u \in (0, \eta), v>0\}.$$ Then we claim that for any $x \in X_{\eta}$, the set $\s_{\eta}(x)$ of slopes of vectors of $\La_x$ in $V_{\eta}$ defined in (\ref{eq:slopeset}) corresponds precisely to the set of times $$\{s >0: h_s x \in \Om_{\eta}\},$$ and that this set is non-empty.  This claim follows from the observation that if $(u,v)^T \in V_{\eta} \cap \La_x$, then, setting $s = \frac{v}{u}$ to be the slope of the vector, we have $$h_{s} (u,v)^T = (u, 0)^T,$$ with $0 < u \le \eta$, that is, $h_s x \in \Om_{\eta}$. This shows that the set of slopes $\s_{\eta}(x)$ is \emph{contained} in the set of return times $\{s >0: h_s x \in \Om_{\eta}\}.$ For the reverse containment, suppose $h_s x \in \Om_{\eta}$. Then there is a vector $(u, 0)^T \in \La_{h_s x}$ with $0 < u \le \eta$, so $h_{-s} (u, 0)^T = (u, su) \in \La_{x} \cap V_{\eta}$, so $s \in \s_{\eta}(x)$ as desired. To show that this set is discrete, note that if there was an accumulation point, say $s_0$, there would be a sequence of distinct vectors in $\La_x \cap V_{\eta}$ accumulating along the line of slope $s_0$, so by passing to a subsequence, we would obtain a sequence of distinct vectors in $\La_x \cap V_{\eta}$ converging to a point along the line segment $\{(u, s_0 u): 0 < u \le \eta\}$, a contradiction to the discreteness of $\La_x$. 

Finally, to show that this set is non-empty, we need to use the Minkowski and vertical strip properties. We first show that for any $x \in X_{\eta}$, there is an $s \in \R$ so that $h_s x \in \Om_{\eta}$.To show this, let $$A = [-\eta, \eta] \times \left[-\frac{c}{4\eta}, \frac{c}{4\eta}\right],$$ where $c = c(x)$ is as in the Minkowski property. Then $A$ is convex, centrally symmetric, and has volume $c$. Since we assume $x \in X_{\eta}$, $\La_x$ does not contain a vertical vector of length at most $\frac{c}{4\eta}$, so we have a vector $(v_1, v_2)^T \in A \cap \La_x$, with $v_1 \neq 0$. Since $\La_x$ is centrally symmetric, we can assume $v_1 >0$, and setting $s = \frac{v_2}{v_1}$, we have $(v_1, 0) \in h_s x$, so $h_s x \in \Om_{\eta}$ as desired. 

Next, we will use the vertical strip property to show that for any $\om \in \Om_{\eta}$, there is an $s'>0$ so that $h_{s'} \om \in \Om_{\eta}$. For $\om \in \Om_{\eta}$, we can use the vertical strip property to find $(u, v)^T \in \La_{\om} \cap V_{\eta}$, and setting $s' = \frac{v}{u}$, we have our result. This completes the proof of Lemma~\ref{lemma:omt}.\qed

\subsubsection{Self-Similarity}\label{sec:selfsim} We let $R_{\eta} : \Om_{\eta} \rightarrow \R^+$ be the return time function $$R_{\eta}(\om) = \min\{s>0: h_s \om \in \Om_{\eta}\},$$ and $T_{\eta}: \Om_{\eta} \rightarrow \Om_{\eta}$ be the return map $$T_{\eta}(\om) = h_{R_{\eta}(\om)} \om.$$ Given two positive parameters $\eta_1, \eta_2 >0$, the maps $T_{\eta_1}$ and $T_{\eta_2}$ are related by the action of the diagonal subgroup $g_t$, via the conjugation relation $g_t h_s g_{-t} = h_{se^{-t}}.$ Precisely, letting $T = 2 \log \frac{\eta_1}{\eta_2},$ we have that $g_{T} \Om_{\eta_2} = \Om_{\eta_1}$, since the horizontal vectors of length (at most) $\eta_2$ get expanded by $e^{T/2} = \frac{\eta_1}{\eta_2}$ to vectors of length (at most) $\eta_1$, and the slopes of the vectors get scaled by $e^T = \left(\frac{\eta_1}{\eta_2}\right)^2.$ This yields the self-similarity relations $$ T_{\eta_1}(g_T\om) = g_T T_{\eta_2}(\om), $$ and $$R_{\eta_1}(g_T \om) = \left(\frac{\eta_1}{\eta_2}\right)^2 R_{\eta_2}(\om)$$ for any $\om \in \Om_{\eta_2}$.

\subsubsection{Generic points} We now prove Theorem~\ref{theorem:meta}. Recall notation: $\sigma_{x, S}$ denotes the Lebesgue probability measure $\frac{ds}{S}$ on the orbit $\{h_s x: 0 \le s \le S\}$. Our assumption is that $$\sigma_{x, S} \rightarrow \nu$$ as $S \rightarrow \infty$, where convergence is in the weak-$\star$ topology (in fact, all of our convergence of measures will be in this topology, so we do not mention it again), and that $\nu$ gives measure $0$ to the set of vertical $x$. We want to show (\ref{eq:slopelim}), that is, for any $0 \le a \le b \le \infty$, $$\lim_{N \rightarrow \infty} \frac{| \G_N^{\eta}(x) \cap [a, b]|}{N} = \tilde{\nu}\left( R_{\eta}^{-1}([a,b])\right),$$ where $$\G_N^{\eta}(x)=\{1 \le i \le N: s_{i+1}^{\eta}(x) - s_i^{\eta}(x)\}$$ is the set of the first $N$ gaps of $\s_{\eta}(x)$, and $\tilde{\nu}$ is the measure on $\Om_{\eta}$ such that $d\nu = d\tilde{\nu}ds$. Let $\om_0= h_{s_1^{\eta}(x)} \in \Om_{\eta}$ be the first hitting point of $\Om_{\eta}$ for the orbit $\{h_s x\}_{s >0}$. Then the assumption $\sigma_{x, S} \rightarrow \nu$ implies that the measures $\sigma_{\om_0, N}$ supported on the given by $$\sigma_{\om_0, N} := \frac{1}{N} \sum_{i=0}^{N-1} \delta_{T_{\eta}^i (\om_0)}$$ satisfy $$\sigma_{\om_0, N} \rightarrow \tilde{\nu}$$ as $N \rightarrow \infty$. Let $\chi_{a, b}$ denote the indicator function of the set $R_{\eta}^{-1}([a,b]).$ Then we claim $$\frac{| \G_N^{\eta}(x) \cap [a, b]|}{N} = \sigma_{\om_0, N} (\chi_{a,b}),$$ since $$R_{\eta}(T^{i} (\om) ) = s_{i+1}^{\eta}(x) - s_{i}^{\eta}(x),$$ which follows from the proof of Lemma~\ref{lemma:omt}. That is, 

\begin{obs} The gap distribution is a Birkhoff sum for the roof function $R_{\eta}$ along the orbit of the map $T_{\eta}$\end{obs}

\noindent Then, applying the convergence $\sigma_{\om_0, N} \rightarrow \tilde{\nu}$ to the bounded function $\chi_{a,b}$, we have our result. The above argument also shows how Corollary~\ref{cor:meta:ergodic} follows from the Birkhoff ergodic theorem applied to the system $(\Om_{\eta}, T_{\eta}, \tilde{\nu})$ and the function $\chi_{a,b}$.\qed\medskip

\subsubsection{Periodic points} To prove Theorem~\ref{theorem:meta:period}, we first need to prove Lemma~\ref{lemma:period}, that is, we need to show if $x$ is $h_s$-periodic, then the set of slope gaps eventually stabilizes. Let $s_0$ be the (minimal) period of $x$, so $h_{s_0} x = x$. Then $h_{s_0} \La_x = \La_x$. Let $$M =\min\{n>0: s_{n}^{\eta}(x) \geq s_0 \}.$$ Then $$\s_{\eta}(h_{s_0}x) = \{ s_M^{\eta}(x) - s_0 < s_{M+1}^{\eta}(x) -s_0 < s_{M+2}^{\eta}(x) -s_0< \ldots < s_{M+i}^{\eta}(x) -s_0< \ldots \},$$ but $h_{s_0}x = x$, so $$s^{\eta}_{M+i}(x) - s_0 = s_{i+1}^{\eta}(h_{s_0} x) = s_{i+1}^{\eta}(x),$$ and thus $$s_{i+1}^{\eta}(x) - s_{i}^{\eta}(x) = s_{M+i}^{\eta}(x) - s_{M+i-1}^{\eta}(x),$$ that is, the gaps have period $M-1$. Setting $N_0 = M-1$, we have, as desired, for any $N \geq N_0$, $$\G_N^{\eta}(x) = \G_{N_0}^{\eta}(x).$$ This completes the proof of Lemma~\ref{lemma:period}.  To complete the proof of Theorem~\ref{theorem:meta:period}, we note that fixing an arbitrary $\eta_0 >0$, periodic orbits for $h_s$ naturally correspond to periodic orbits for the return map $T_{\eta_0}$. The above argument shows  the periodic orbit of $T_{\eta_0}$ associated to the orbit $$\{ h_s g_{-T} x: 0 \le s \le e^T s_0\}$$ has length $N_0(\eta)$, where $T = 2\log(\frac{\eta}{\eta_0}).$ By moving $x$ by $h_s$, and applying Lemma~\ref{lemma:omt}, we can assume $x \in \Om_{\eta_0}$. We can then define the probability measure $\sigma_{x, \eta}$ supported on the periodic orbit for $T_{\eta_0}$  on $\Om_{\eta_0}$ induced by the probability measure $\rho_{x, T} = \frac{ds}{s_0 e^T}$ supported on the orbit $\{ h_s g_{-T} x: 0 \le s \le e^T s_0\}$, via $$\sigma_{x, \eta} := \frac{1}{N_0(\eta)} \sum_{i=0}^{N_0(\eta) -1} \delta_{T_{\eta_0}^i (x)}.$$ Recall from \S\ref{sec:selfsim} that $g_T$ scales $V_{\eta}$ to $V_{\eta_0}$ and multiplies slopes (and gaps of slopes) of vectors by the factor $\left(\frac{\eta}{\eta_0}\right)^2$. Thus we have $$\frac{\left|\left(\frac{\eta}{\eta_0}\right)^2 \G_{N_0(\eta)}^{\eta}(x) \cap [a, b]\right|}{N_0(\eta)} = \sigma_{x, \eta} (\chi_{a,b}),$$ and applying our convergence assumption, $\rho_{x, T} \rightarrow \nu$, which implies $$\sigma_{x, \eta} \rightarrow \tilde{\nu},$$ we have our result.\qed

\subsection{Circle Limits}\label{subsec:circle:meta} The proof of Theorem~\ref{theorem:meta:circle} is essentially already contained in~\cite[\S9]{MS} and in~\cite[\S4]{AChaika}, which we discussed in \S\ref{sec:affine} and \S\ref{sec:saddle}. We sketch the argument once again: we define, for $x \in X$, $\sigma>0$, $$p_{i}(x, \sigma) : = \lim_{R \rightarrow \infty} \lambda(\theta: |A_{\theta}^{\sigma}(R) \cap \La_x| = i).$$ By the renormalization argument described in \S\ref{sec:affine}, we have, for $R>>0$, $$|\Lambda_{x} \cap A_{\theta}^{\sigma}(R)| \approx |g_{t} r_{-\theta} \Lambda_{x} \cap T(\sigma)|,$$ where $t = 2\log R$. Thus, $$p_i(x, \sigma) \approx \lambda_{x, R}( \{y \in X: |\Lambda_{y} \cap T(\sigma)| = i).$$ Applying our convergence assumption $\lambda_{x, R} \rightarrow \nu$, we have our result.\qed

\section{Further Questions}\label{sec:quest} \noindent We collect some further questions and speculations. We discuss the space of translation surfaces in \S\ref{sec:trans:ques}, the space of affine lattices in \S\ref{sec:affine:ques}, and speculate wildly about other problems in \S\ref{sec:other:ques}.

\subsection{Translation surfaces}\label{sec:trans:ques} It would be intriguing to push the machinery developed in this note further than the set of lattice surfaces. In particular, it would be interesting to check which $SL(2,\R)$ (or $h_s$)-invariant measures on $\Om_g$ are supported on the set of surfaces satisfying the vertical strip condition, and then to attempt to explicitly identify the transversal and the return time function (and the associated distribution). It would be particularly nice to do this for the Masur-Veech measure on $\Om_g$, and to understand if there were any `continuity' properties of the gap distributions- namely, given a sequence of surfaces $\om_n \rightarrow \om$, do the gap distributions of saddle connection directions converge in any appropriate sense? Also, are there any surfaces $\om \in \Om_g$ which do have `truly random' behavior, that is, an exponential distribution of gaps? We conjecture that this is not possible. 

\subsection{Affine Lattices}\label{sec:affine:ques} A nice test case, which in fact corresponds to an $SL(2, \R)$-invariant subset of $\Om_2$ (see \cite[\S6]{AChaika}), is the collection of \emph{affine unimodular lattices}, as discussed in \S\ref{sec:affine}. If we could give an explicit description of the (4-dimensional) transversal $\Om_{\eta} \subset \tilde{X_2} = \left(SL(2, \R) \ltimes \R^2\right)/\left(SL(2, \Z) \ltimes \Z^2\right)$, and understand the return time function, we would be able to recover the (2-dimensional) results of~\cite{MS}, explicitly using unipotent flows.

\subsection{Other gap distributions}\label{sec:other:ques} As discussed in the introduction, although our paper was inspired by~\cite{Elkies}, our machinery does not seem to immediately give the gap distribution for the sequence $\{ \{\sqrt n\} \}_{n \geq 0}$. It would be a nice application to formulate the results of loc. cit. in the language our machine. A more ambitious, but probably very difficult project would be to try and understand the apparent exponential distribution of the gaps  in $\{ \{n^{\alpha}\} \}_{n \geq 0}$ for other $0 < \alpha < 1$. Another sequence of great interest in applications (see~\cite[\S9]{Hodgson}) would be the gaps in the sequence of squares of lengths of lattice vectors $\{|a+b\tau|^2: a, b \in \Z\}$ for any fixed $\tau \in \C$, which are also conjectured to be exponential. In fact, we would be very excited to see any application of our technology which lead to a non-exotic (i.e., exponential) distribution of gaps.

\end{document}